\documentclass[graybox]{svmult}
\usepackage{amsmath,amsfonts,amssymb,mathrsfs} %,eufrak}
\usepackage{graphicx}
\newcommand*{\abs}[1]{\left\vert #1\right\vert}
\newcommand*{\dif}{\mathrm{d}}

\newcommand*{\e}{\mathrm{e}}
\newcommand*{\mi}{\mathrm{i}}
\newcommand*{\sca}[1]{\left\langle #1 \right\rangle }
\newcommand*{\cnj}[1]{\overline{#1}}
\newcommand*{\hardy}[1]{\mathsf{H}^{#1}}
\newcommand*{\hardyT}[1]{\mathsf{H}^{#1}({\mathbb T})}
\newcommand*{\hardyR}[1]{\mathsf{H}^{#1}({\mathbb R})}
\newcommand*{\uhp}{{\mathbb H}}
\newcommand*{\blp}{{\mathsf B}}
\newcommand*{\mat}[1]{{\mathsf #1}}
\begin{document}

\title*{Multiscale decompositions of Hardy spaces}
\author{Ronald R.~Coifman and Jacques Peyri\`ere}
\institute{Ronald R.~Coifman \at Department of Mathematics, Program in
  Applied Mathematics, Yale University, New Haven, CT 06510, USA,
  \email{coifman-ronald@yale.edu} \and Jacques Peyri\`ere \at Institut
  de Math\'ematiques d'Orsay, CNRS, Universit\'e Paris-Saclay, 91405
  Orsay, France, \email{jacques.peyriere@universite-paris-saclay.fr}}

% \keywords{Blaschke factorization, Phase unwinding, Takenaka basis,
%   Hardy spaces, inner function, invariant subspaces}

%\subjclass[2000]{30B50; 30A10, 42C40, 65T99 }

\maketitle
\section{Introduction}
We would like to elaborate on a program of analysis pursued by Alex Grossmann and his collaborators on the analytic utilization of the phase of Hardy functions, as a multiscale signal processing tool. 

An inspiration at the origin of "wavelet" analysis  (when Grossmann, Morlet, Meyer and collaborators  were interacting and exploring  versions of multiscale representations) was provided, by the  analysis of  holomorphic   signals, for which, the images of the phase of Cauchy wavelets were remarkable in their ability to reveal intricate singularities or dynamic structures, such as instantaneous frequency jumps, in musical recordings.  
This work which was pursued by Grossmann, Kronland Martinet et al \cite{sound}  exploiting phase and amplitude variability of holomorphic signals was challenged by computational complexity as well as by the lack of simple, efficient, mathematical processing, and generalizations to higher dimensional signals. It was mostly bypassed by the orthogonal wavelet transforms. We aim to show that these ideas are powerful nonlinear subtle tools.

  Our goal here is to follow their seminal work and introduce recent developments in nonlinear analysis. In particular we will sketch methods extending conventional Fourier analysis, exploiting both phase and amplitudes of holomorphic functions.

  The miracles of nonlinear complex analysis, such as factorization and composition of functions lead to new versions of holomorphic wavelets, and  relate them to multiscale dynamical systems.
  
  Our story  interlaces the role of the phase of signals with their analytic/geometric properties.  The Blaschke factors are a key ingredient, in building analytic tools, starting with the Malmquist Takenaka orthonormal bases of the Hardy space $\hardyT{2}$, continuing with "best" adapted bases obtained through phase unwinding, and concluding with relations to composition of Blaschke products and their dynamics (on the disc, and on invariant subspaces of $\hardyT{2}$).
  Specifically we discuss multiscale orthonormal holomorphic wavelet bases, related to  Grossmann's and Morlet's program~\cite{wavelet}, and associated   generalized scaled holomorphic orthogonal bases, to dynamical systems, obtained by composing Blaschke factors.  
  
   We also, remark, that the phase of a Blaschke product is a one layer neural net with ($\arctan$ as an activation sigmoid) and that the composition is a "Deep Neural Net" whose depth is the number of compositions, our results provide a wealth of related libraries of orthogonal bases .
  
  We  sketch these ideas in various "vignette" subsections and
 refer for more details on analytic methods~\cite{CP}, related to the Blaschke based nonlinear phase 
unwinding decompositions~\cite{coifman,CSW,nahon}, we also consider
orthogonal decompositions of invariant subspaces of Hardy spaces. In
particular we constructed a multiscale decomposition, described below, of the Hardy space of the upper half-plane. 

Such a decomposition can
be carried in the unit disk by conformal mapping. A somewhat different
multiscale decomposition of the space $\hardyT{2}$ has been
constructed by using Malmquist-Takenaka bases associated with Blaschke
products whose z\'eros are
$\displaystyle (1-2^{-n})\e^{2\mi\pi j/2^n}$ where $n\ge 1$ and
$0\le j< 2^n$ \cite{feichtinger}. Here we provide a variety  of multiscale decompositions by considering iterations of Blaschke products.

\section{Preliminaries and notation}

For $p\ge 1$, $\hardyT{p}$ stands for the space of analytic functions~$f$
on the unit disk~${\mathbb D}$ such that
\begin{equation*}
\sup_{0< r<1} \int_{0}^{2\pi} |f(r\e^{\mi\theta})|^p\frac{\dif
  \theta}{2\pi} < +\infty.
\end{equation*}
Such functions have boundary values almost everywhere, and the Hardy
space $\hardyT{p}$ can be identified with the set of $L^p$ functions
on the torus~${\mathbb T}=\partial{\mathbb D}$ whose Fourier
coefficients of negative order vanish.

A subspace of~$\hardyT{p}$ is \emph{invariant} if it is invariant under 
multiplication by~$\e^{\mi\theta}$ (or by~$z$), depending whether these
functions are considered as functions on~${\mathbb T}$ or~${\mathbb
  D}$. An inner function is a bounded analytic function on the unit
disk whose boundary values have modulus~1 almost everywhere. It is
known that the invariant subspaces are of the form $u\hardyT{p}$
where~$u$ is an inner function~\cite{helson,hoffman}. The inner function~$u$ is determined
by the invariant subspace up to multiplication by a constant of
modulus~1.

If $f$ and $g$ are two functions on ${\mathbb T}$ (in $L^p$ and
$L^{p/(p-1)}$ for some $p\in [1,+\infty)$), let
\begin{equation*}
\sca{f,g} = \frac{1}{2\pi} \int_0^{2\pi}
f(\e^{\mi\theta})\cnj{g(\e^{\mi\theta})} \,\dif \theta.
\end{equation*}

Let $\mathscr H$ be the operator of orthogonal projection of $L^2({\mathbb T})$
onto $\hardyT2$. 
It results from the properties of the Hilbert transform that this operator
extends as a bounded operator from $L^p({\mathbb T})$ to $\hardyT{p}$
for $1< p< +\infty$.

If $u$ is an inner function, let $\chi_u$ be the operator of
multiplication by~$u$ (which is an isometry of all the $L^P$). Then the
operator ${\mathscr H}_u=\chi_u{\mathscr H}\chi_u^{-1}$ is the operator of orthogonal projection
of $L^2$ onto $u\hardyT2$. It results that this operator extends as a
bounded operator from $L^p({\mathbb T})$ to $\hardyT{p}$ for all $p\in
(1,+\infty)$ with a norm independent of~$u$. In other terms, for all~$p>
1$, there exists~$C_p$ such that, for all~$u$ and all $f\in
L^{p}({\mathbb T})$,
\begin{equation}\label{projection}
\|{\mathscr H}_uf\|_p \le C_p\|f\|_p.
\end{equation}
\bigskip

There is a parallel theory for analytic functions on the upper half
plane $\uhp = \{x+\mi y \ :\ y>0\}$. The space of analytic
functions~$f$ on $\uhp$ such that
$$\sup_{y>0} \|f(\cdot+\mi y)\|_{L^p({\mathbb R})} < +\infty$$
is denoted by $\hardyR{p}$. These functions have boundary values in
$L^p({\mathbb R})$ when $p\ge 1$. The space $\hardyR{p}$ is identified
to the space of $L^p$ functions whose Fourier transform vanishes on
the negative half line~$(-\infty,0)$.

A subspace of $\hardyR{2}$ is said to be invariant if it is stable by
multiplication by the functions $\e^{2\mi\pi\xi x}$ for all
$\xi>0$. As previously, the invariant subspaces are of the form
$u\,\hardy{2}$ where $u$ is an inner function, i.e., a bounded analytic
function on~$\uhp$ whose boundary values are of modulus~1 almost
everywhere.

As previously, the operators of orthogonal projections on invariant
subspaces extend, for any $p\in (1,+\infty)$, as continuous operators
on $\hardyR{p}$ with a uniform bound for their norms.

\section{Malmquist-Takenaka bases on the torus}

\begin{lemma}\label{takenaka}
  Let $a$ be a complex number of modulus less than~1 and $u$ be an
  inner function. Then $(z-a)u\hardy{2}$ has codimension~1 in
  $u\hardy{2}$ and $\displaystyle \frac{\sqrt{1-|a|^2}}{1-\cnj{a}z}u$
  is a unit vector in the orthogonal complement of $(z-a)u\hardy{2}$
  in $u\hardy{2}$.
\end{lemma}

\proof Since $f\mapsto uf$ is an isometry of $\hardy{2}$ onto
$u\hardy{2}$, it is enough to consider the case~$u=1$. One has
\begin{eqnarray*}
  \sca{(z-a)f(z),\frac{\sqrt{1-|a|^2}}{1-\cnj{a}z}} &=&
  \frac{1}{2\pi}
  \int_{-\pi}^{\pi} (\e^{\mi\theta}-a)f(\e^{\mi\theta})
  \frac{\sqrt{1-|a|^2}}{1-a\e^{-\mi\theta}}\,\dif \theta\\
&=& \frac{\sqrt{1-|a|^2}}{2\pi}\int_{-\pi}^{\pi}
\e^{\mi\theta}f(\e^{\mi\theta})\,\dif \theta = 0. 
\end{eqnarray*}

Also, if $f$ is orthogonal to $(1-\cnj{a}z)^{-1}$ one has
\begin{equation*}
0 = \frac{1}{2\pi} \int_{-\pi}^{\pi}
\frac{f(\e^{\mi\theta})}{1-a\e^{-\mi\theta}}\,\dif \theta\\
= \frac{1}{2\mi\pi} \oint \frac{f(z)}{z-a}\,\dif z = f(a),
\end{equation*}
so $f\in (z-a)\hardy{2}$.\medskip

Now let $(a_n)_{n>0}$ be a sequence of complex numbers of modulus less
than~1.
For $n\ge 0$, let
 $$\blp_n(z) = \prod_{0\le j< n} \frac{z-a_j}{1-\cnj{a}_jz} \text{\quad
  and\quad} \phi_n(z) = \blp_n(z)\frac{\sqrt{1-|a_n|^2}}{1-\cnj{a}_nz}.$$
It results from Lemma~\ref{takenaka} that, if
$$\sum_{n\ge 1}(1-|a_j|^2) = +\infty,$$
the functions $\phi_n$ form an orthonormal basis of $\hardy{2}$.

If $\displaystyle\sum_{n\ge 1}(1-|a_j|^2) < +\infty$ the functions
$\phi_n$ form an orthonormal basis of $\hardy{2}\ominus \blp\hardy{2}$,
where $\blp$ is the convergent Blaschke product
$$\blp(z) = \prod_{j>0} \frac{\cnj{a}_j}{|a_j|}\frac{z-a_j}{1-\cnj{a}_jz}.$$

Consider a sequence $(\blp_m)_{m\ge 1}$:
$$\blp_m(z) = \prod_{j>0}
\frac{\cnj{a}_{m,j}}{|a_{m,j}|}\frac{z-a_{m,j}}{1-\cnj{a}_{m,j}z}$$ of
convergent Blaschke products such that
$\displaystyle \sum_{m,j>0} (1-|a_{m,j}|^2) = +\infty$.

Let $\blp_0=1$ and, for $n\ge 0$ and $m\ge 1$,
\begin{eqnarray*}
\blp_{m,n}(z) &=& \prod_{0\le j< n} \frac{z-a_{m,j}}{1-\cnj{a}_{m,j}z}\\
  \phi_{m,n}(z) &=&
         \blp_{m-1}(z)\blp_{m,n}(z)\frac{\sqrt{1-|a_{m,n}|^2}}{1-\cnj{a}_{m,n}z}.
\end{eqnarray*}
Then, $(\phi_{m,n})_{n\ge 1}$ is an orthonormal basis of
$\blp_{m-1}\hardy{2}\ominus \blp_{m-1}\blp_m\hardy{2}$, and
$(\phi_{m,n})_{m\ge 1, n\ge 1}$ is an orthonormal basis of
$\hardy{2}$.

The bases so obtained are the \emph{Malmquist-Takenaka bases}~\cite{takenaka}.

% Theorem~\ref{convergence} implies that, if $1< p<+\infty$
% and $f\in \hardy{p}$, the series
% $\displaystyle \sum_{n\ge 0} \sca{f,\phi_n}\phi_n$ converges to~$f$ in
% $\hardy{p}$.

\section{The upper half plane}

We present some prior results~\cite{CP}, without proof.
 In this section one simply writes $\hardy{2}$ instead of
$\hardyR{2}$.

\subsection{Malmquist-Takenaka bases}

Let
$(a_j)_{1\le j}$ be a sequence (finite or not)) of complex numbers
with positive imaginary parts and such that
\begin{equation}\label{upper} 
\sum_{j\ge 0} \frac{\Im a_j}{1+|a_j|^2} < +\infty.
\end{equation}
The corresponding Blaschke product is
$$\blp(x) = \prod_{j\ge 0}
\frac{\abs{1+a_j^2}}{1+a_j^2}\,\frac{x-a_j}{x-\cnj{a}_j},$$
where, $0/0$, which appears if $a_j=\mi$, should be understood as~1.
The factors $\displaystyle \frac{\abs{1+a_j^2}}{1+a_j^2}$ insure the
convergence of this product when there are infinitely many
zeroes. But, in some situations, it is more convenient to use other
convergence factors as we shall see below.

Whether the series~\eqref{upper} is convergent or not, one defines
(for $n\ge 0$) the functions
\begin{equation*}
\phi_n(x) = \frac{1}{\sqrt{\pi}}\left( \prod_{0\le j< n}
\frac{x-a_j}{x-\cnj{a}_j}\right)\, \frac{1}{x-\cnj{a}_n}.
\end{equation*}
Then these functions form a orthonormal system in $\hardy{2}$. If the
series~\eqref{upper} diverges, it is a basis of $\hardy{2}$, otherwise
it is a basis of the orthogonal complement of $\blp\,\hardy{2}$ in
$\hardy{2}$.

\subsection{A multiscale Wavelet decomposition}

The infinite products
\begin{equation*}
G_n(x) = \prod_{j\le n} \frac{j-\mi}{j+\mi}\, \frac{x-j-\mi}{x-j+\mi}
\text{\quad and\quad } G(x) = \prod_{j\in {\mathbb Z}} \frac{j-\mi}{j+\mi}\,
\frac{x-j-\mi}{x-j+\mi}
\end{equation*}
 can be expressed in terms of known functions:
\begin{equation*}
G_n(x) = \frac{\Gamma(-\mi-n)}{\Gamma(\mi-n)}\,
\frac{\Gamma(x-n+\mi)}{\Gamma(x-n-\mi)} \text{\quad and\quad } G(x) =
\frac{\sin \pi(\mi-x)}{\sin \pi(\mi+x)}.
\end{equation*}

\subsection{An orthonormal system}

Consider the function
$$\phi(x) = \frac{\Gamma(x-1+\mi)}{\sqrt{\pi}\Gamma(x-\mi)}.$$
It is easily checked that
$$\phi(x-n) =
\frac{\Gamma(\mi-n)}{\Gamma(-\mi-n)} \, \frac{G_n(x)}{\sqrt{\pi}\bigl(
  x-(n+1)+\mi\bigr)}.$$
Set $\phi_n(x) = \phi(x-n)$. For fixed~$m$, the functions
$\phi_n/G_m$, for $n\ge m$, form a Malmquist-Takenaka basis of
$(G/G_m)\hardy{2}$. In other terms, the functions $\phi_n$, for $n\ge
m$, form an orthonormal basis of $G_m\hardy{2}\ominus
G\hardy{2}$. This means that the functions $\phi_n$ (for $n\in
{\mathbb Z}$) form a Malmquist-Takenaka basis of the orthogonal
complement of $G\hardy{2}$ in $\hardy{2}$.

\subsubsection{Multiscale decomposition}

As $|1-G(2^nx)|\le C2^{n}$ all the products
\begin{equation}
{\mathscr B}_n(x) = \prod_{j< n} G(2^jx)
\end{equation}
are convergent and $\displaystyle \lim_{n\to -\infty} {\mathscr B}_n =
1$ uniformly.

Let ${\mathscr B}={\mathscr B}_0$. Obviously, ${\mathscr B}_n(x) =
{\mathscr B}(2^nx)$.

Consider the following subspaces of $\hardy{2}$:
\begin{equation*}
{\mathsf E}_n = {\mathscr B}_n\hardy{2}.
\end{equation*}
This is a decreasing sequence. The space $\displaystyle {\mathsf
  E}_{+\infty} = \bigcap_{n\in {\mathbb Z}} {\mathsf E}_n$ is equal
to~$\{0\}$ since a function orthogonal to this space would have too
many zeros, and the space $\displaystyle {\mathsf E}_{-\infty} =
\mathrm{closure~of}\bigcup_{n\in {\mathbb Z}} {\mathsf E}_n$ is
equal to $\hardy{2}$ since ${\mathscr B}_n$ converges uniformly to 1
when $n$ goes to $-\infty$.

For all $n$ and~$j$, let
$$\phi_{n,j}(x) =
2^{n/2}\phi(2^nx-j){\mathscr B}(2^{n}x).$$
Then, for all~$n$, $(\phi_{n,j})_{j\in {\mathbb Z}}$ is an orthonormal
basis of ${\mathsf E}_{n}\ominus {\mathsf E}_{n+1}$. At last
$(\phi_{n,j})_{n,j\in{\mathbb Z}}$ is an orthonormal basis of
$\hardy{2}$.

\section{ Adapted MT bases, "phase unwinding"}
Our goal is to find  a "best"  adapted Malmquist Takenaka basis to analyze a given function the idea is to peel off the oscillation of a function by dividing by its  Blaschke product, this procedure is iterated to yield an expansion in an orthogonal collection of  functions or Blaschke products which of course are naturally embedded in a MT basis, once the zeroes are ordered.

\subsection{The unwinding series.} There is a natural way to  iterate the Blaschke factorization, it is inspired by the power series expansion
of a holomorphic function. If $G$ has no zeroes inside $\mathbb{D}$, its Blaschke factorization is the trivial one $G = 1 \cdot G$, however,
the function $G(z)-G(0)$ certainly has at least one root inside the unit disk $\mathbb{D}$ and will therefore yield some nontrivial Blaschke factorization
$G(z) - G(0) = \blp_1 G_1$.
We write
\begin{align*}
 F &= \blp \cdot G\\
&= \blp \cdot (G(0) + (G(z) - G(0))) \\
&= \blp \cdot (G(0) + \blp_1 G_1) \\
&= G(0) \blp + \blp \blp_1 G_1.
\end{align*}
An iterative application gives rise to  the \textit{unwinding series}
$$ F = a_1 \blp_1 + a_2 \blp_1 \blp_2 + a_3 \blp_1 \blp_2 \blp_3 + a_4 \blp_1 \blp_2 \blp_3 \blp_4+ \dots$$
This orthogonal expansion first appeared in the PhD thesis of Michel Nahon~\cite{nahon} and independently by T.~Qian in \cite{qtao,qw}. Given a general function $F$ it is not numerically feasible to actually compute the roots of the function; a crucial insight in \cite{nahon} is that this is not necessary -- one can numerically obtain the Blaschke product
in a stable way by using a method that was first mentioned in a paper of Guido and Mary Weiss \cite{ww} and has been investigated with respect to stability by Nahon \cite{nahon} Using the boundedness of the Hilbert transform one can prove easily convergence in $L^p, 1<p<\infty$ .
 
\subsection { The fast algorithm of  Guido and Mary Weiss \cite{ww}  }

Our starting point is the theorem that any Hardy function   can be
decomposed as
$$ F = \blp \cdot G,$$
where $B$ is a Blaschke product, that is a function of the form
$$ \blp(z) = z^m\prod_{i \in I}{\frac{\overline{a_i}}{|a_i|}\frac{z-a_i}{1-\overline{a_i}z}},$$
where $m \in \mathbb{N}_{0}$ and $a_1, a_2, \dots \in \mathbb{D}$ are zeroes inside the unit disk $\mathbb{D}$ 
and $G$ has no roots in $\mathbb{D}$. For $|z|=1$ we have $|\blp(z)| = 1$
which motivates the analogy 
$$\blp \sim \mbox{frequency and}~G \sim \mbox{amplitude}$$
for the function restricted to the boundary. However, the function $G$ need not be constant: it can be any function that never vanishes
inside the unit disk. If $F$ has roots inside the unit disk, then the Blaschke factorization $F = \blp\cdot G$ is going to be nontrivial (meaning
$\blp \not\equiv 1$ and $G \not\equiv F$). $G$ should be 'simpler' than $F$ because the winding number around the origin decreases.

In fact since $|F|=|G|$ and $\ln(G)$ is analytic in the disk we have formally that

$G=\exp(\ln|F|+ \mi(\ln|F|)^\sim))=\exp(\mathscr H(\ln|F|))$   

and   $\blp=F/G$ 

G can be computed easily using the FFT \cite{nahon}.

\section{Iteration of Blaschke products.}
%\subsubsection*{A few introductory remarks.}

% If~$f$ is a function, whenever it makes sense, $\iter{f}{n}$ stands
% for the $n$th iterate of~$f$: $\iter{f}{0} = \Id$, $\iter{f}{n+1} =
% \iter{f}{n}\circ f$.\bigskip

 % If $u_1$ and $u_2$ are inner functions on the unit disk ${\mathbb D}$ so
 % is $u_2\circ u_1$.
 We are interested in the case we iterate finite Blaschke products:
\begin{equation*}
\blp(z) = \e^{\mi\theta}z^\nu \prod_{j=1}^\mu \frac{z+a_j}{1+\cnj{a}_jz}, 
\end{equation*}
where~$\mu$ and~$\nu$ are nonnegative integers and the $a_j$ are complex
numbers of modulus less than~1.

It is well known that ${\mathbb T}$ and ${\mathbb D}$ are globally
invariant under~$\blp$, as well as the complement of~$\cnj{\mathbb D}$ in
the Riemann sphere.

We have
\begin{equation*}
\e^{-\mi\varphi}\blp(\e^{\mi\varphi}z) =
\e^{\mi\bigl(\theta+(\nu+\mu-1)\varphi\bigr)} z^\nu \prod_{j=1}^\mu
\frac{z+b_j}{1+\cnj{b}_jz},
\end{equation*}
where $b_j=\e^{-\mi\varphi} a_j$. This means that for iteration
purpose we may assume~$\theta=0$.

When~$\nu\ge 1$, it results from the Schwarz lemma that~$0$ is the
unique fixed point of~$\blp$ in~${\mathbb D}$ and that this fixed point
is attracting. The basin of attraction of~$0$ contains a disk centered
at~$0$. The sequence of iterates $\blp_n$ is a normal family, since it converges
towards~0 on some neighborhood of~$0$, it has a unique limit
point. This means that~$\blp_n$ converges to~$0$ uniformly on any
compact subset of~${\mathbb D}$.

Now, suppose~$\nu=0$. We have
\begin{equation*}
\blp(z)\cnj{\blp(1/\cnj{z})}=1,
\end{equation*}
so, if ~$\alpha$ is a fixed point of~$\blp$, so is $1/\cnj{\alpha}$.

we have two possibilities:
\begin{enumerate}
\item
  There exists a fixed point $\alpha\in {\mathbb D}$ of~$\blp$. Then
  if $\phi(z) = \frac{z+\alpha}{1+\cnj{\alpha}z}$, then~$0$ is a fixed
  point of the Blaschke product~$\phi^{-1}\circ \blp\circ\phi$. It
  results from the preceding discussion that~$\alpha$ is the unique
  fixed point of~$\blp$ in~${\mathbb D}$, that it is attracting, and that
  the sequence of iterates of~$\blp$ converges to the
  constant~$\alpha$.
\item All the fixed points lie in~${\mathbb T}$. Then according to the
  Denjoy-Wolff theorem~\cite{carleson}, again the sequence of iterates
  of~$\blp$ converges to some constant.
\end{enumerate}

 Moreover, the fixed points of~$\blp$ are the roots of an equation of
 degree~$\mu+1$, so if $\blp$ has a fixed point in~${\mathbb D}$, it has
 $\mu-1$ fixed points in~$\mathbb T$ (none of them being attracting). In this
 case, the dynamics on $\mathbb T$ is that of a cookie cutter.

 When  $\blp_1$ and $\blp_2$ are finite Blaschke products with $n_1$ and
 $n_2$ zeros, then $\blp_2\circ \blp_1$ is a finite Blaschke products with
 $n_1n_2$ zeros (of course zeros are counted with their
 multiplicities). This is obvious by considering the variation of
 arguments on the boundary of the disk.

 Let $\blp$ be a finite Blaschke product having at least two
 zeros. Then, one may consider the dynamical system which it
 generates. Let $\blp_{0}(z)=z$, and $\blp_{n+1} =
 \blp_{n}\circ\blp$.

 \begin{lemma}\label{diverge}
 Consider a Blaschke product $F$ of the form $F(z) = z\blp(z)$, where
 $\blp$ is nonconstant. Then there exists a sequence
 $0,a_1,a_2,\dots,a_j,\dots$ of complex numbers in the unit disk and
 an increasing sequence $(\nu_j)_{j\ge1}$ of positive integers such
 that $a_1,a_2,\dots,a_{\nu_n}$ are the zeros, counted according to
 their multiplicity, of~$F_{n}$. Moreover $\displaystyle \sum_{j\ge
   1}(1-|a_j|) = +\infty$.
\end{lemma}

\proof A simple recursion shows
\begin{equation}\label{dvl}
   F_{n} = z(\blp\!\circ\!\! F)(\blp\!\circ\!\! F_2)\cdots(\blp\!\circ\!\! F_{n-1}).
\end{equation}
This proves the existence of the sequences $(a_j)$ and $(\nu_j)$. In
fact, if $\blp$ has~$k-1$ zeros, $\nu_n = k^n$.

If $\blp(0)=0$, then the multiplicity of 0 in $F_n$ increases
with~$n$, so $$\sum_{j\ge 1}(1-|a_j|) = +\infty.$$

Otherwise, let us compute the derivative $F_{n}'(0)$. On the one hand,
due to~\eqref{dvl} and $F(0)=0$, we have
$F_n'(0)=\blp(0)^n$. Therefore $\displaystyle \lim_{n\to \infty}
F_n'(0)=0.$ On the other hand, $$F_n'(0) = \lim_{z\to 0}
\frac{F_n(z)}{z}, \text{\quad so\quad} \abs{F_n'(0)} =
\prod_{j=1}^{\nu_n} |a_j|.$$
Thus $\displaystyle \prod_{j\ge 1}|a_j| = 0$. This ends the proof.

\begin{corollary}\label{Diverge}
Let $F$ be a finite Blaschke product with a fixed point~$\alpha$
inside the unit disk. Then there exits a sequence
$\alpha,a_1,a_2,\dots,a_j,\dots$ of complex numbers in the unit disk
and an increasing sequence $(\nu_j)_{j\ge1}$ of positive integers such
that $a_1,a_2,\dots,a_{\nu_n}$ are the zeros, counted according to
their multiplicity, of~$F_{n}$. Moreover $\displaystyle \sum_{j\ge
  1}(1-|a_j|) = +\infty$.
\end{corollary}

\proof Just consider $\varphi_\alpha\circ F\circ\varphi_\alpha$, with
$\displaystyle \varphi_\alpha(z) = \frac{\alpha-z}{1-\cnj{\alpha}z}$,
and use Lemma~\ref{diverge}.

When $F=z\blp$, if $|w|<1$ and $F(z)=w$, we have $|w|= |z||\blp(z)|<
|z|$. A refinement of this observation leads in some circumstances to
an estimate of the speed of accumulation of the zeros to the
boundary. Here is an example.

\begin{lemma}\label{bounds} Let $k$ be a positive integer and
  $a\in {\mathbb D}$, $a\ne 0$. Let $\displaystyle F(z) =
  z\frac{z^k-a^k}{1-\cnj{a}^kz^k}$. Then there exist two increasing
  and concave functions~$g$ and~$h$ on $[0,1]$ such that
  $g(1)=h(1)=1$, $g(0)=0$, $g(t)<h(t)$ if $t<1$, and such
  that, if $F(z)=w$ with $|w|\ge |a|$, we have $g(|w|)< |z|< h(|w|)$.
\end{lemma}

\proof One can check the following equalities, valid for $0<\rho, r<
1$ and any real~$\varphi$:
$$\frac{r+\rho}{1+r\rho} = \inf_{\theta}
\abs{\frac{r\e^{\mi\theta}-\rho\e^{\mi\varphi}}{1-r\rho\e^{\mi(\theta-\varphi)}}}
<\sup_{\theta}
\abs{\frac{r\e^{\mi\theta}-\rho\e^{\mi\varphi}}{1-r\rho\e^{\mi(\theta-\varphi)}}}
= \frac{\abs{r-\rho}}{1-r\rho}.$$

Let $\rho=|a|$. Suppose $F(z)=w$, with $|w|=t\ge \rho$, and $|z| =
r$. As already observed $|z|> |w|$. Then the preceding equalities give
$$r\frac{r^k+\rho^k}{1+r^k\rho^k}\le t\le r
\frac{r^k-\rho^k}{1-r^k\rho^k}.$$

Let $\varphi(r) = r\frac{r^k-\rho^k}{1-r^k\rho^k}$ and
$\psi(r) = r\frac{r^k+\rho^k}{1+r^k\rho^k}$. It is routine to check
that $\varphi$ and $\psi$ is increasing and convex, $\psi$ from
$[0,1]$ onto~$[0,1]$ and $\varphi$ from $[\rho,1]$ onto $[0,1]$. Then
the inverse functions~$g$ and~$h$ of~$\psi$ and~$\phi$ will do.

\begin{proposition}\label{sandwich} Keep the notation and hypotheses
  of Lemma~\ref{bounds}. Then, for all~$w$ such that $F_{n+1}
  (w)/F_n(0)=0$, we have $g_n(|a|)<|w|<h_n(|a|)$.
\end{proposition}

\proof The zeros of $F_{n+1}/F_n$ are inverse image by~F of the zeros
of $F_n/F_{n-1}$. Also the zeros of $F/F_0$ are $\e^{2\mi
  j\pi/k}a$. Knowing that, if $z\ne 0$ we have $|F(z)|< |z|$, we
conclude that all the zeros (except~0!) have a modulus larger than or
equal to~$|a|$.

Lemma~\ref{bounds} shows that for all zero $w$ of $F_2/F_1$ we have
$g(|a|)\le |w|\le h(|a|)$. Then we can perform a recursion. Assuming
that all zeros of $F_n/F_{n-1}$ have an absolute value between
$g_{n-1}(|a|)$ and $h_{n-1}(|a|)$, Lemma~\ref{bounds} gives that the
absolute value of the zeros of $F_{n+1}/F_n$ all lie between
$g_n(|a|)$ and $h_n(|a|)$.

It is worth noticing that~1 is an attracting fixed point of both~$g$
and~$h$. Besides we have
$$g'(1) = \frac{1}{1+k\frac{1-|a|^k}{1+|a^k|}} \text{\quad and \quad
} h'(1) = \frac{1}{1+k\frac{1+|a|^k}{1-|a^k|}}.$$

\begin{figure}[h]
\begin{center}
\includegraphics[width=38mm]{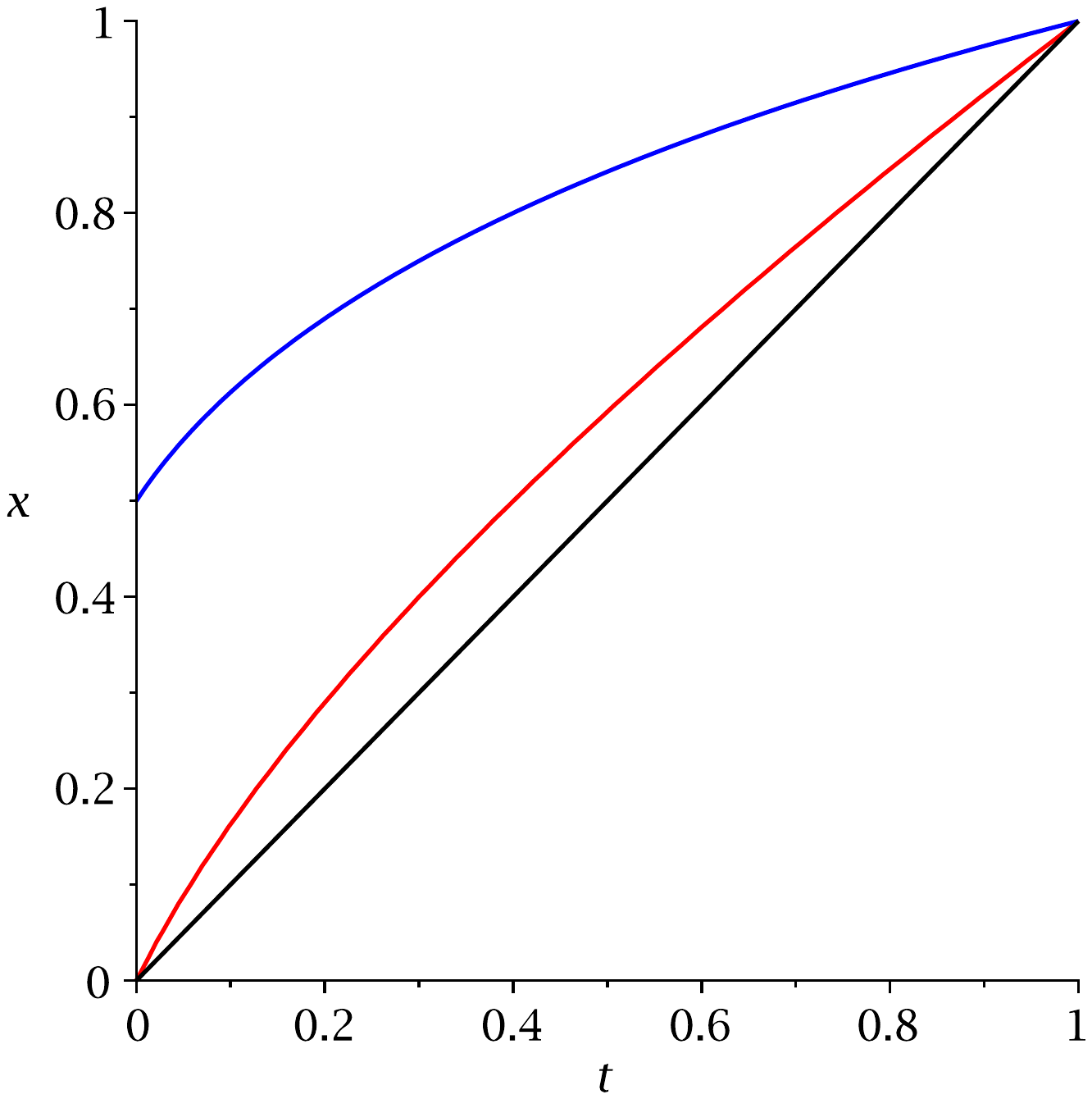}
\includegraphics[width=38mm]{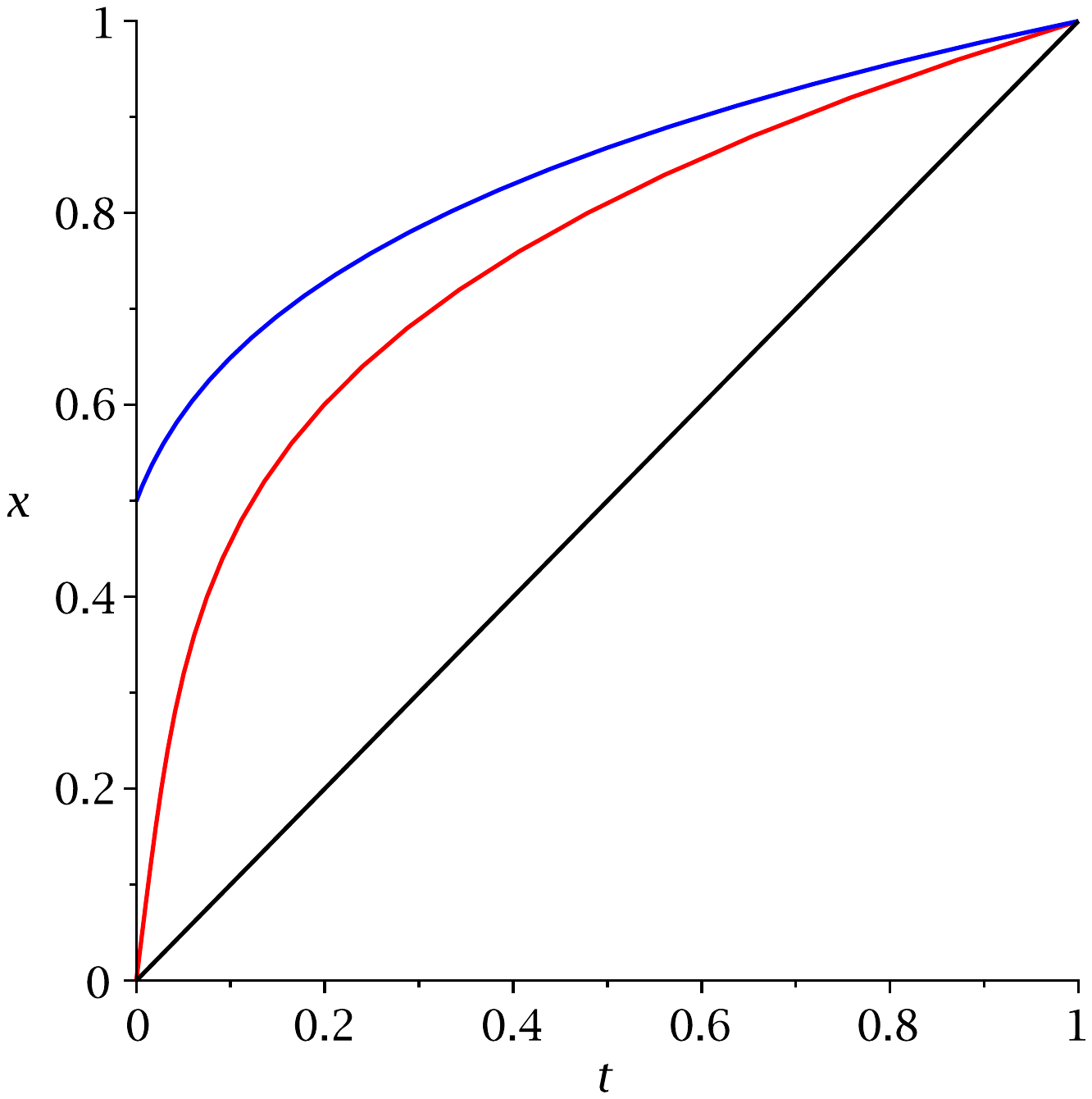}
\includegraphics[width=38mm]{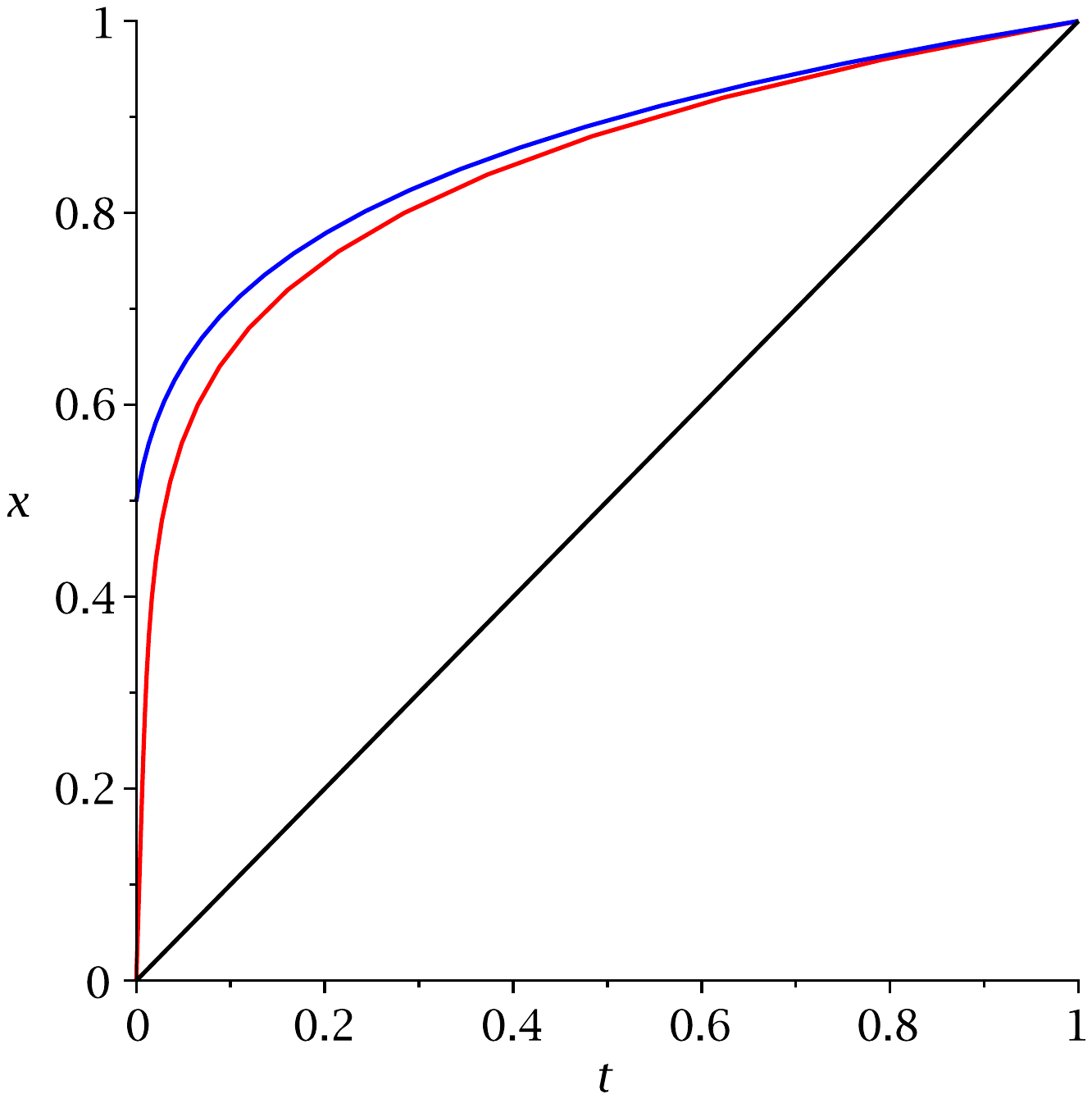}
\end{center}
\caption{The functions $h$ and $g$, with $a=0.5$ and $k=1,3,5$}
\end{figure}
 
\subsection{Multiscale decomposition}

Each Blaschke product $\blp$ defines invariant subspaces of
$\mat{H}^p$. The projection on this space is given by the kernel
$\displaystyle \frac{\blp(z)\cnj{\blp(w)}}{z-w}$. This projection is
continuous for $1< p< +\infty$.

Let $F$ be a Blaschke product of degree at least~2 with a fixed point
inside the unit disk. Its iterates define a hierarchy of nested
invariant subspaces ${\mathsf E}_n = F_n\mat{H}^2$.

Due to Corollary~\ref{Diverge},
$\displaystyle\bigcap_{n\ge 1} {\mathsf E}_n = \{0\}$.

The Takenaka construction provides orthonormal bases of $E_n\ominus
E_{n+1}$. But this is not canonical as it depends on an ordering of the
zeros of $F_{n+1}/F_n$.

Figure~2 shows 1st, 3rd, and 5th iterates of $F(z)
=z(z-2^{-1})/(1-2^{-1}z)$. Figure~3 displays the same for $F(z)
=z(z-2^{-2})/(1-2^{-2}z)$. The upper pictures display the phases
modulo $2\pi$ (values in the interval $(-\pi,\pi]$) of theses Blaschke
products while the lower pictures display minus the logarithms of
their absolute value. The coordinates $(x,y)$ correspond to the point
$\e^{-y + \mi x}$. On these figures it is easy to locate the zeros,
specially by looking at the phase which then has an abrupt jump.

\begin{figure}[h]
\begin{center}\hspace*{-20mm}
  \includegraphics[width=130mm,height=60mm]{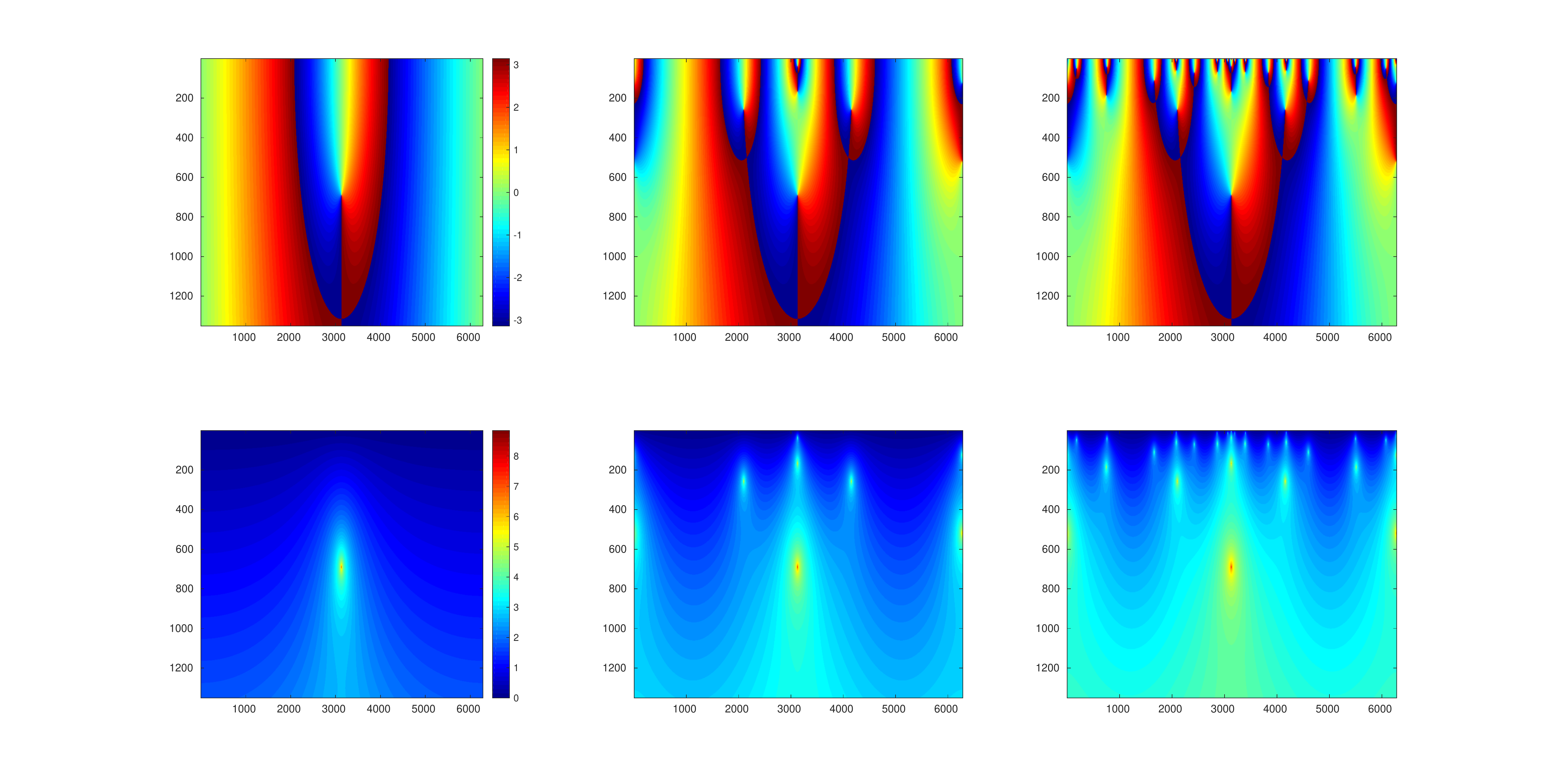}\\
\caption{ The argument and the absolute value of $F(z)$, $F^{(3)}(z)$,
  and $F^{(5}(z)$, with $F(z)=\frac{z(z-2^{-1})}{1-2^{-1}z}$ and $z=
  \exp(-y+\mi x)$.}
\end{center}\end{figure}

\begin{figure}[h]
\begin{center}
\hspace*{-20mm}
\includegraphics[width=140mm,height=60mm]{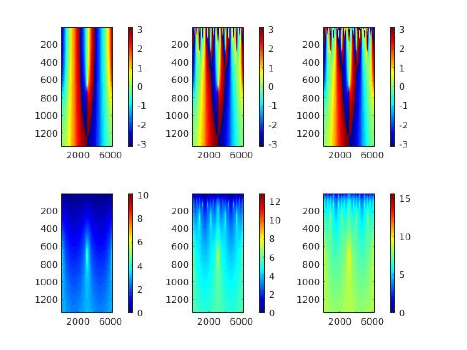}
\caption{ The argument and the absolute value of $F(z)$, $F^{(3)}(z)$,
  and $F^{(5)}(z)$, with $F(z)=\frac{z(z^2-2^{-2})}{1-2^{-2}z^2}$ and $z=
  \exp(-y+\mi x)$.}
\end{center}
\end{figure}

\section{Remarks on Iteration of Blaschke products as a "Deep Neural Net"}
 In the upper half plane let
$(a_j)_{1\le j}$ be a finite sequence  of complex numbers
with positive imaginary parts.

The corresponding Blaschke product on the line is
$$\blp(x) = \prod_{j\ge 0}
(x-a_j)/(x-\cnj{a}_j),$$

 We can write $\blp(x)=\exp(\mi\theta(x))$. where 
$$\theta(x)=\sum_{j\ge 0}\sigma[(x-\alpha_j)/(\beta_j)] $$
where $a_j=\alpha_j+\mi\beta_j$  and $\sigma=\arctan(x) + \pi/2 $ is a sigmoid.

This is precisely the form of a single layer in a Neural Net, each unit has a weight and bias determine by $a_j$ . We obtain the various layers of a deep net through the composition 
of each layer with a preceding layer. In our preceding examples we took a single short layer given  by a simple  Blaschke term with two zeroes in the first layer that we iterated to obtain an orthonormal Malmquist Takanaka basis ( we could have composed  different elementary  products at each layer), demonstrating the versatility of the method to generate highly complex functional representations.

As an example let $F(z)$ be mapped from G  (4.2) in the section on wavelet construction.
$$ F(z)=G(w)= \frac{\sin(\pi(\mi-w))}{\sin(\pi(\mi+w))}$$
with $w=\frac{\mi(1-z)}{(1+z)}$  

we can view the phase of F as a neural layer which when composed with itself results in a Phase which is a two layer neural net represented graphically in fig 4. 

Where each end of a color droplet corresponds to one zero or unit of the two layer net. \begin{figure}[h]
\begin{center}
\hspace*{-10mm}
\includegraphics[width=120mm,height=120mm]{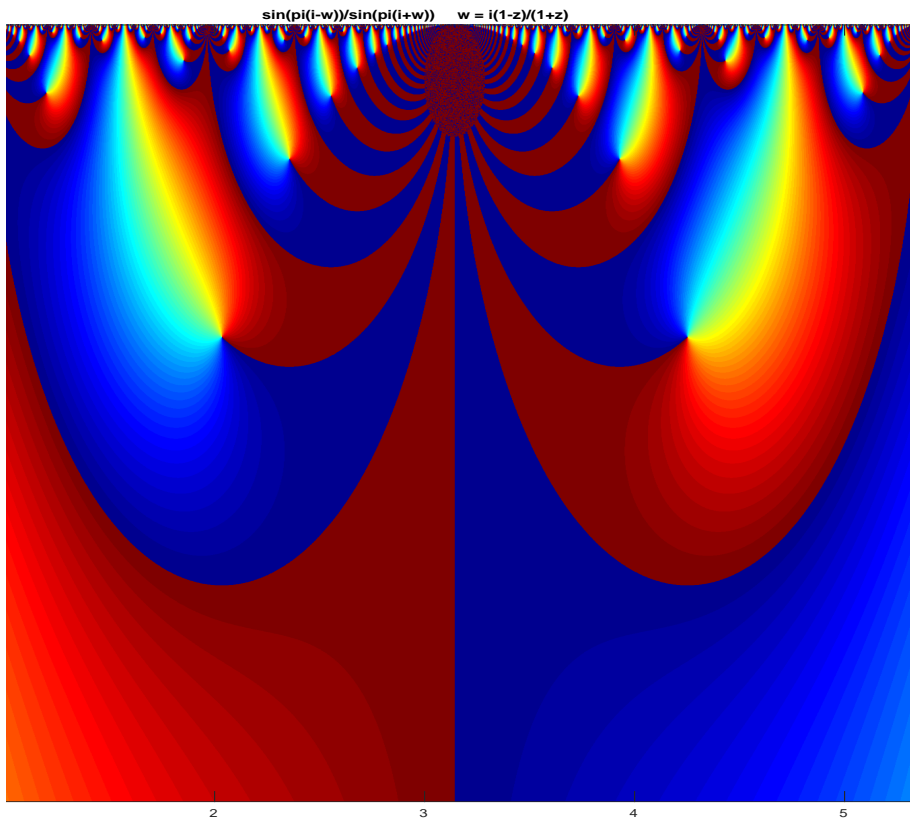}
\caption{ Two iterations  of $F(z)$}
\end{center}
\end{figure}

We conclude by referring to Daubechies et al \cite{ReluDNN}for a description of a similar iteration for piecewise affine functions in which simple affine functions play the role of a Blaschle product. 

\section*{ Appendix:An example}
In this section $a\in {\mathbb D}$ and $\displaystyle \blp(z) =
\left(\frac{z+a}{1+\cnj{a}z}\right)^2$.
\begin{proposition}
\begin{enumerate}
\item
  If $27\abs{a}^4-18\abs{a}^2+ 8\Re{a} < 1$, then $\blp$ has a unique
  fixed point inside~${\mathbb D}$ and a unique fixed point
  in~${\mathbb T}$. Moreover $\blp_n$ converges, uniformly on any
  compact subset of~${\mathbb D}$, towards the attracting fixed point.
\item If $27\abs{a}^4-18\abs{a}^2+ 8\Re{a} \ge 1$, $\blp$ has no fixed
  point in~${\mathbb D}$; it has three fixed points in~${\mathbb T}$
  and only one of them is attracting.  Moreover $\blp_n$
  converges, uniformly on any compact subset of~${\mathbb D}$, towards
  the attracting fixed point.
\end{enumerate}
\end{proposition}

\proof The fixed points of $B$ are the roots of the equation
\begin{equation}\label{fp}
  \cnj{a}^2z^3+(2\cnj{a}-1)z^2-(2a-1)z-a^2=0.
\end{equation}

If we write $z=x+\mi y$ and $a=t+\mi u$, with $x,\,y,\,t$, and $u$
real, this equation becomes

\begin{multline*}
  - t^2x^3 + 3t^2xy^2 - 6tux^2y + 2tuy^3 + u^2x^3 - 3u^2xy^2 - 2tx^2 +
  2ty^2\\ - 4uxy + t^2 + 2tx - u^2 - 2uy + x^2 - y^2 -
  x\\ +\mi(-3t^2x^2y + t^2y^3 + 2tux^3 - 6tuxy^2 + 3u^2x^2y - u^2y^3 -
  4txy\\ + 2ux^2 - 2uy^2 + 2tu + 2ty + 2ux + 2xy - y) = 0.
\end{multline*}

So~$x$ and~$y$ are the real solutions to a systen of two polynomial
equations. If we look for the fixed points in~${\mathbb T}$ we have to
add the equation $x^2+y^2=1$. By substituting $1-x^2$ to $y^2$ in the
real and imaginary parts of the previous equation we get the following
two equations
\begin{multline}\label{eq1}
(2x - 1)(2t^2x - 2u^2x + t^2 - u^2 + 2t - 1)y\\
= - 2u(x + 1)(2x - 1)(2tx - t + 1)
\end{multline}
\noindent and
\begin{multline}\label{eq2}
  2u(2x + 1)(2tx - t + 1)y\\
  = + (2x + 1)(x - 1)(2t^2x - 2u^2x + t^2 - u^2 + 2t - 1).
\end{multline}

By eliminating~$y$ between Equations~\eqref{eq1} and~\eqref{eq2} we
get the polynomial
\begin{multline*}
P = {4(t^2 + u^2)^2}x^3 + (8t^3 + 8tu^2 - 4t^2 + 4u^2)x^2\\ + (-3t^4 -
6t^2u^2 - 3u^4 - 4t^3 + 12tu^2 + 6t^2 + 2u^2 - 4t + 1)x\\ - (t^2 + 2tu
- u^2 + 2t - 2u - 1)(t^2 - 2tu - u^2 + 2t + 2u - 1).
\end{multline*}

If $x_0$ is a real root of~$P$, then $y_0$ computed by
substituting~$x_0$ to~$x$ in Equation~\eqref{eq1} or in
Equation~\eqref{eq2} is real and satisfies $x_0^2+y_0^2=1$. So the
real roots of~$P$ have absolute value less than or equal to~1.

The real roots of~$P$ are the real parts of the fixed points of~$\blp$ on
the circle~${\mathbb T}$. So, the number of fixed points of modulus~1,
depends on the sign of the discriminant of $P$ viewed as a
polynomial in~$x$. This discriminant is
$$ -(t^2 + u^2 - 1)^2\bigl(27(t^2+u^2)^2-18(t^2+u^2)+8t-1\bigr).
$$
So all depends on the sign of $Q = 27(t^2+u^2)^2-18(t^2+u^2)+8t-1$,
i.e., the sign of $27|a|^4-18|a|^2+8\Re{a}-1$. This accounts for the
first assertion and part of the second one.\bigskip

Now we give an algebraic proof that, according to the Denjoy-Wolff theorem, when~$Q>0$ there is one fixed point where
$|B'|<1$.

%% We have $\displaystyle 1-wB'(z) = R/(\mi uz - tz - 1)^3$ with
%% \begin{multline*}
%% R = w({\mi}u - t)^3z^3 - 3w({\mi}u - t)^2z^2 + (3{\mi}uw - 2t^2 - 3tw
%% - 2u^2 + 2)z\\ - 2{\mi}t^2u - 2{\mi}u^3 - 2t^3 - 2tu^2 + 2{\mi}u + 2t
%% - w.
%% \end{multline*}

Let $P_1$ be the polynomial in~$z$ obtained by replacing~$a$ by $t+\mi
u$ in the left-hand side of Equation~\eqref{fp}. The resultant
of~$P_1$ and~$1-wB'(z)$ considered as polynomials in~$z$, is $(\mi u - t)^3(t^2 + u^2 - 1)^4R(w)$, where
\begin{multline*}
R(w) = 8(t^2 + u^2)(t^2 + u^2 - 1)w^3 + (12t^4 + 24t^2u^2 + 12u^4 - 4t^2 - 4u^2 + 8t)w^2\\ + 2(3t^2 + 3u^2 + 1)(t^2 + u^2 - 1)w + (t^2 + u^2 - 1)^2.
\end{multline*}

The roots of~$R$, considered as a polynomial in~$w$, are the inverses of
the derivative of~$\blp$ evaluated at the fixed points of~$\blp$.

Let $R_1(w) = R(1+w)$. We have
\begin{multline*}
R_1(w) = 8(t^2 + u^2)(t^2 + u^2 - 1)w^3 \\+ 4\bigl(
9(t^2+u^2)^2-7(t^2+u^2)+2t\bigr)w^2+2Qw+Q.
\end{multline*}

We see that the coefficients of~$R_1$ present one variation of
sign. This means that this polynomial has a positive root. Therefore,
due to Descartes' rule, $R$ has a root larger than one. This means
that there is a fixed point where the derivative of~$\blp$ is in the
interval~$(0,1)$.

The following figure shows the regions corresponding to the previous
discussion.  
\begin{figure}[h]
  \begin{center}
    \caption{ When $(t,u)$ lies between the
      cardioid and the circle, all the fixed points have modulus~1, when it
      lies inside, there is a fixed point inside~${\mathbb D}$.}
    
    \includegraphics[width=90mm]{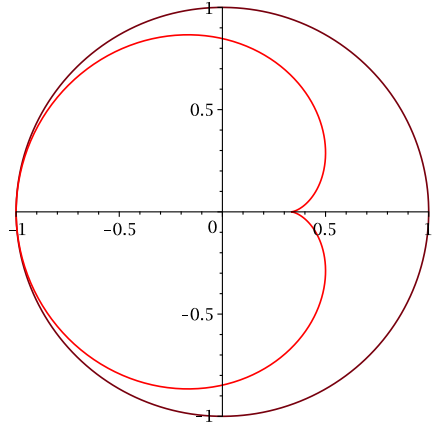}
  \end{center}
\end{figure}

The equation of the red curve (a cardioid) is  $
27(t^2+u^2)^2-18(t^2+u^2)+8t-1=0$.\medskip

When the attracting fixed point is on the boundary of the disk, $|\blp_n(0)|$ converges exponentially fast towards~1. Therefore, if $(z_j)_{j\ge 0}$ is the sequence of the zeroes (counted according to their multiplicities) of all the iterates of~$\blp$, one has $$\sum_{j\ge 0} (1-|z_j|) < +\infty. $$

\end{document}